\documentclass[12pt]{article}
\usepackage{a4}
\usepackage{amsfonts}
\usepackage{amssymb}

\begin{document}

\title{\textbf{Lower bound for the poles of Igusa's p-adic zeta functions}}
\author{Dirk Segers}
\date{August 30, 2005}

\maketitle

\begin{abstract}
Let $K$ be a $p$-adic field, $R$ the valuation ring of $K$, $P$
the maximal ideal of $R$ and $q$ the cardinality of the residue
field $R/P$. Let $f$ be a polynomial over $R$ in $n>1$ variables
and let $\chi$ be a character of $R^{\times}$. Let $M_i(u)$ be the
number of solutions of $f=u$ in $(R/P^i)^n$ for $i \in
\mathbb{Z}_{\geq 0}$ and $u \in R/P^i$. These numbers are related
with Igusa's $p$-adic zeta function $Z_{f,\chi}(s)$ of $f$. We
explain the connection between the $M_i(u)$ and the smallest real
part of a pole of $Z_{f,\chi}(s)$. We also prove that $M_i(u)$ is
divisible by $q^{\ulcorner(n/2)(i-1)\urcorner}$, where the corners
indicate that we have to round up. This will imply our main
result: $Z_{f,\chi}(s)$ has no poles with real part less than
$-n/2$. We will also consider arbitrary $K$-analytic functions
$f$.
\end{abstract}

\section{Introduction}
\noindent \textbf{(1.1)} Let $K$ be a $p$-adic field, i.e., an
extension of $\mathbb{Q}_p$ of finite degree. Let $R$ be the
valuation ring of $K$, $P$ the maximal ideal of $R$, $\pi$ a fixed
uniformizing parameter for $R$ and $q$ the cardinality of the
residue field $R/P$. For $z \in K$, let $\mathrm{ord} \, z \in
\mathbb{Z} \cup \{+\infty\}$ denote the valuation of $z$,
$|z|=q^{-\mathrm{ord} \, z}$ the absolute value of $z$ and
$\mbox{ac} \, z =z \pi^{-\mathrm{ord} \, z}$ the angular component
of $z$.

Let $\chi$ be a character of $R^{\times}$, i.e., a homomorphism
$\chi: R^{\times} \rightarrow \mathbb{C}^{\times}$ with finite
image. We formally put $\chi(0)=0$. Let $e$ be the conductor of
$\chi$, i.e., the smallest $a \in \mathbb{Z}_{>0}$ such that
$\chi$ is trivial on $1+P^a$.

\vspace{0,5cm}

\noindent \textbf{(1.2)} Let $f$ be a $K$-analytic function on an
open and compact subset $X$ of $K^n$ and put $x=(x_1,\ldots,x_n)$.
Igusa's $p$-adic zeta function of $f$ and $\chi$ is defined by
\[
Z_{f,\chi}(s)=\int_X \chi(\mathrm{ac} \, f(x)) |f(x)|^s \, |dx|
\]

\vspace{0,5cm} \footnoterule{ \noindent \footnotesize{2000
\emph{Mathematics Subject Classification}. Primary 11D79 11S80;
Secondary 14B05}
\\ \emph{Key words.} Polynomial congruences, Igusa's $p$-adic zeta
function.}
\newpage

\noindent for $s \in \mathbb{C}$, $\mbox{Re}(s) \geq 0$, where
$|dx|$ denotes the Haar measure on $K^n$, so normalized that $R^n$
has measure $1$. Igusa proved that it is a rational function of
$q^{-s}$, so that it extends to a meromorphic function
$Z_{f,\chi}(s)$ on $\mathbb{C}$ which is also called Igusa's
$p$-adic zeta function of $f$. We will write $Z_{f,\chi}(t)$ if we
consider $Z_{f,\chi}(s)$ as a function in the variable $t=q^{-s}$.
If $\chi$ is the trivial character, we will also write $Z_f(s)$
and $Z_f(t)$.

\vspace{0,5cm}

\noindent \textbf{(1.3)} A power series $f=\sum_{(i_1,\ldots,i_n)
\in \mathbb{Z}_{\geq 0}^n} c_{i_1,\ldots,i_n}x_1^{i_1} \ldots
x_n^{i_n}$ over $K$ is convergent in $(a_1,\ldots,a_n) \in K^n$ if
and only if $|c_{i_1,\ldots,i_n}a_1^{i_1} \ldots a_n^{i_n}|
\rightarrow 0$ if $i_1+\cdots+i_n \rightarrow \infty$. If $f$ is
convergent at every $(a_1,\ldots,a_n) \in (P^k)^n$ for some $k \in
\mathbb{Z}$, then $f$ is called a convergent power series.

Because a $K$-analytic function is locally described by convergent
power series, we only have to consider this type of functions in
the study of Igusa's $p$-adic zeta function. By performing a
dilatation of the form $(x_1,\ldots,x_n) \mapsto
(\pi^kx_1,\ldots,\pi^kx_n)$, we may moreover suppose that $f$ is
rigid, i.e., convergent on $R^n$. The coefficients of a rigid
$K$-analytic function $f$ on $R^n$ have the property that
$|c_{i_1,\ldots,i_n}| \rightarrow 0$ if $i_1+\cdots+i_n
\rightarrow \infty$. Consequently, $|c_{i_1,\ldots,i_n}|$ is
bounded and we can multiply $f$ by an element of $K$ to obtain a
series over $R$. So we only have to study $Z_f(s)$ for rigid
$K$-analytic functions $f$ on $R^n$ defined over $R$. See also
\cite[Chapter 2]{Igusabook}.

\vspace{0,5cm}

\noindent \textbf{(1.4)} Let $f$ be a rigid $K$-analytic function
on $R^n$ defined over $R$. Igusa's $p$-adic zeta function of such
an $f$ has an important connection with congruences. For $i \in
\mathbb{Z}_{\geq 0}$ and $u \in R/P^i$, let $M_i(u)$ be the number
of solutions of $f(x) \equiv u \mbox{ mod } P^i$ in $(R/P^i)^n$.
Put $M_i:=M_i(0)$.

The $M_{i+e}(\pi^iu)$, $u \in (R/P^e)^{\times}$, describe
$Z_{f,\chi}(t)$ through the relation
\[
Z_{f,\chi}(t) = \sum_{i=0}^{\infty} \sum_{u \in (R/P^e)^{\times}}
\chi(u) M_{i+e}(\pi^iu) q^{-n(i+e)}t^i.
\]

If $\chi$ is the trivial character, all the $M_i$'s describe and
are described by $Z_f(t)$ through the relation
\[ Z_f(t) = P(t) - \frac{P(t)-1}{t}, \]
where the Poincar\'e series $P(t)$ of $f$ is defined by
\[ P(t) = \sum_{i=0}^{\infty} M_i(q^{-n}t)^i. \]
Remark that $P(t)$ is a rational function because $Z_f(t)$ has
this property.

\vspace{0,5cm}

\noindent \textbf{(1.5)} Igusa's $p$-adic zeta function is often
studied by using an embedded resolution of $f$. The well known
fact that $Z_{f,\chi}(s)$ has no poles with real part less than
$-1$ if $n=2$ is easily proved in this way. We used this method in
\cite{Segers} to determine all the values less than $-1/2$ which
occur as the real part of a pole of some $Z_{f,\chi}(s)$ if $n=2$,
and all values less than $-1$ if $n=3$. In particular, we proved
that there are no poles with real part less than $-3/2$ if $n=3$.
In arbitrary dimension $n>1$, we saw in \cite[Section
3.1.4]{Segersthesis} that it is easy to prove that there are no
poles with real part less than $-(n-1)$ and we conjectured that
this bound can be sharpened to $-n/2$.

Let $f$ be a rigid $K$-analytic function on $R^n$ defined over
$R$. In \cite{Segers} we proved that there exists an integer $a$
such that $M_i$ is an integer multiple of $q^{\ulcorner
(n/2)i-a\urcorner}$ for all $i$ if this conjecture is true in
dimension $n$ for the trivial character. Consequently, this
divisibility property of the $M_i$ is true for $n=2$ and $n=3$.
The statement of this property is so easy that we tried to find an
elementary proof, and with success. It generalized easily to
arbitrary dimension and to the more general class of numbers
$M_i(u)$. This is the subject of the second section. We deduce
there that $M_i(u)$ is divisible by
$q^{\ulcorner(n/2)(i-1)\urcorner}$ for all $i \in
\mathbb{Z}_{>0}$.

\vspace{0,5cm}

\noindent \textbf{(1.6)} The poles of Igusa's $p$-adic zeta
function are an interesting object of study for example because
they are related to the monodromy conjecture
\cite[(2.3.2)]{Denefreport}. In the third section, we explain the
connection between the $M_i(u)$ and the smallest real part of a
pole of Igusa's $p$-adic zeta function. Let $l$ be the smallest
real part of a pole of $Z_f(s)$. We proved in \cite{Segers} that
there exists an integer $a$ which is independent of $i$ such that
$M_i$ is an integer multiple of $q^{\ulcorner(n+l)i-a\urcorner}$
for all $i \in \mathbb{Z}_{\geq 0}$. We repeat this proof for
completeness and we also prove the converse: if there exists an
integer $a$ such that $M_i$ is an integer multiple of
$q^{\ulcorner (n+l')i-a \urcorner}$ for all $i \in
\mathbb{Z}_{\geq 0}$, then $l' \leq l$. The last statement has an
analogue if we are dealing with a character. Together with (1.5),
this will imply that $Z_{f,\chi}(s)$ has no pole with real part
less than $-n/2$. Remark that this bound is optimal: $Z_f(s)$ has
a pole in $-n/2$ if $f$ is equal to
$x_1x_2+x_3x_4+\cdots+x_{n-1}x_n$ for $n$ even and
$x_1x_2+x_3x_4+\cdots+x_{n-2}x_{n-1}+x_n^2$ for $n$ odd, see
\cite[Corollary 10.2.1]{Igusabook}.

\section{A theorem on the number of solutions of congruences}
\noindent \textbf{(2.1)} Let $f$ be a rigid $K$-analytic function
on $R^n$ defined over $R$. Let $i \in \mathbb{Z}_{\geq 0}$ and $u
\in R/P^i$. The numbers $M_i(u)$ and $M_i$, which are defined in
(1.4), will sometimes be denoted by respectively $M_i(f,u)$ and
$M_i(f)$.


\vspace{0,5cm}

\noindent \textbf{(2.2)} Let $f$ be a rigid $K$-analytic function
on $R^n$ defined over $R$. Let $(b_1,\ldots,b_n) \in R^n$. Then
$g(y_1,\ldots,y_n):=f(b_1+y_1,\ldots,b_n+y_n)$ is a rigid
$K$-analytic function on $R^n$ defined over $R$. Consequently,
$h(z_1,\ldots,z_n):=g(\pi z_1,\ldots,\pi z_n) = f(b_1+\pi
z_1,\ldots,b_n+\pi z_n)$ is a power series which is convergent on
$\pi^{-1}R \supset R$ and the coefficient of a monomial of degree
$r$ in this power series is in $P^r$.

Note also that the coefficients of a convergent power series are
related with partial derivatives.

\vspace{0,5cm}

\noindent \textbf{(2.3) Theorem.} \textsl{Let $n \in
\mathbb{Z}_{>1}$. Then we have that}
\[
q^{\ulcorner (n/2)(i-1)\urcorner} \mid M_i(f,u)
\]
\textsl{for all rigid $K$-analytic functions $f$ on $R^n$ defined
over $R$, $i \in \mathbb{Z}_{>0}$ and $u \in R/P^i$.}

\vspace{0,2cm}

\noindent \textsl{Remark.} The number
$\ulcorner(n/2)(i-1)\urcorner$ is the smallest integer larger than
or equal to $(n/2)(i-1)$.

\vspace{0,2cm}

\noindent \emph{Proof.} Note that we may suppose that $u$ is zero,
because $f-u$ can be replaced by $f$. So we have to prove that
\[
q^{\ulcorner (n/2)(i-1)\urcorner} \mid M_i(f)
\]
for every rigid $K$-analytic function $f$ on $R^n$ defined over
$R$ and for every $i \in \mathbb{Z}_{>0}$.

The argument is by induction on $i$. For $i=1$, the statement is
trivial. Let $k \in \mathbb{Z}_{\geq 2}$. Suppose that the
statement is true for $i=1,\ldots,k-1$. We prove the statement for
$i=k$. Let $(b_1,\ldots,b_n) \in R^n$. It is enough to prove that
the number of solutions of
\begin{eqnarray}
f(b_1+\pi z_1,\ldots,b_n+\pi z_n) \equiv 0 \mbox{ mod } P^k
\end{eqnarray}
in $(R/P^{k-1})^n$ is a multiple of $q^{\ulcorner
(n/2)(k-1)\urcorner}$. Put $h(z_1,\ldots,z_n) := f(b_1+\pi
z_1,\ldots,b_n+\pi z_n)$. Then $h$ is a rigid $K$-analytic
function on $R^n$ which is defined over $R$. Moreover, the
coefficients of the $z_j$, $j=1,\ldots,n$, are in $P$ and the
coefficients in terms of higher degree are in $P^2$. We explained
this in (2.2).

\textsl{Case 1: Not all the coefficients in the linear part of $h$
are in $P^2$.} Then the number of solutions of (1) in
$(R/P^{k-1})^n$ is equal to $0$ or $q^{(n-1)(k-1)}$. This is
actually Hensel's lemma. Because $(n-1)(k-1) \geq \ulcorner
(n/2)(k-1) \urcorner$, we are done.

\textsl{Case 2: All the coefficients in the linear part of $h$ are
in $P^2$.} Write $h(z_1,\ldots,z_n)$ $=\pi^2
\tilde{h}(z_1,\ldots,z_n)$, where $\tilde{h}$ is a rigid
$K$-analytic function on $R^n$ defined over $R$. Equation (1)
becomes
\begin{eqnarray}
\tilde{h}(z_1,\ldots,z_n) \equiv 0 \mbox{ mod } P^{k-2}.
\end{eqnarray}
We want to prove that the number of solutions of this congruence
in $(R/P^{k-1})^n$ is a multiple of $q^{\ulcorner (n/2)(k-1)
\urcorner}$. If $k=2$, the number of solutions of (2) in $(R/P)^n$
is $q^n$, so we are done because $n \geq \ulcorner n/2 \urcorner$.
If $k>2$, the number of solutions of (2) in $(R/P^{k-1})^n$ is
$q^n M_{k-2}(\tilde{h})$, which is a multiple of $q^n q^{\ulcorner
(n/2)(k-3) \urcorner} = q^{\ulcorner (n/2)(k-1) \urcorner}$. Here
we used the induction hypothesis for $\tilde{h}$ and $i=k-2$.
$\qquad \Box$

\vspace{0,2cm}

\noindent We also give a proof of the theorem which is without
induction.

\vspace{0,2cm}

\noindent \textsl{Alternative proof.} Let $\mathcal{O} \subset
(R/P^i)^n$ be the set of solutions of $f(x_1,\ldots,x_n) \equiv u
\mbox{ mod } P^i$ in $(R/P^i)^n$. We give a partition of
$\mathcal{O}$ with the property that the number of elements of
every subset in this partition is a multiple of $q^{\ulcorner
(n/2)(i-1) \urcorner}$.

Let $r$ be $i/2$ if $i$ is even and $(i-1)/2$ if $i$ is odd. We
associate a subset of $\mathcal{O}$ to every element
$(b_1,\ldots,b_n) \in \mathcal{O}$.

\textsl{Case 1:} $(\partial f / \partial x_j)(b_1,\ldots,b_n)
\equiv 0 \mbox{ mod } P^r$ for every $j \in \{1,\ldots,n\}$. We
associate to $(b_1,\ldots,b_n)$ the set
$\mathcal{O}_{(b_1,\ldots,b_n)} := (b_1,\ldots,b_n)+(P^{i-r})^n$.
The number of elements of $\mathcal{O}_{(b_1,\ldots,b_n)}$ is
$q^{rn}$ and this is a multiple of $q^{\ulcorner (n/2)(i-1)
\urcorner}$ because $rn \geq \ulcorner (n/2)(i-1) \urcorner$.
Moreover, $\mathcal{O}_{(b_1,\ldots,b_n)}$ is a subset of
$\mathcal{O}$ because all the coefficients of $h(z_1,\ldots,z_n):=
f(b_1+\pi^{i-r}z_1,\ldots,b_n+\pi^{i-r}z_n)-u$ are in $P^i$.

\textsl{Case 2:} $(\partial f / \partial x_j)(b_1,\ldots,b_n)
\not\equiv 0 \mbox{ mod } P^r$ \textsl{for at least one} $j \in
\{1,\ldots,n\}$. Let $k$ be the number in $\{0,\ldots,r-1\}$ such
that $(\partial f / \partial x_j)(b_1,\ldots,b_n) \equiv 0 \mbox{
mod } P^k$ for every $j \in \{1,\ldots,n\}$ and such that
$(\partial f /
\partial x_j)(b_1,\ldots,b_n) \not\equiv 0 \mbox{ mod } P^{k+1}$
for some $j \in \{1,\ldots,n\}$. We associate to
$(b_1,\ldots,b_n)$ the subset $\mathcal{O}_{(b_1,\ldots,b_n)} :=
\mathcal{O} \cap ((b_1,\ldots,b_n)+(P^{k+1})^n)$ of $\mathcal{O}$.
The number of elements of $\mathcal{O}_{(b_1,\ldots,b_n)}$ is
equal to the number of solutions of
$f(b_1+\pi^{k+1}z_1,\ldots,b_n+\pi^{k+1}z_n) \equiv u \mbox{ mod }
P^i$ in $(R/P^{i-k-1})^n$. All the coefficients of
$h(z_1,\ldots,z_n):=
f(b_1+\pi^{k+1}z_1,\ldots,b_n+\pi^{k+1}z_n)-u$ are in $P^{2k+1}$,
the coefficient of at least one $z_j$, $j \in \{1,\ldots,n\}$, is
not in $P^{2k+2}$ and the coefficients in terms of degree at least
2 are in $P^{2k+2}$. Write $h(z_1,\ldots,z_n)$ $=\pi^{2k+1}
\tilde{h}(z_1,\ldots,z_n)$. We have to calculate the number of
solutions of $\tilde{h}(z_1,\ldots,z_n) \equiv 0 \mbox{ mod }
P^{i-2k-1}$ in $(R/P^{i-k-1})^n$. This is equal to $q^{nk}
M_{i-2k-1}(\tilde{h})=q^{nk+(n-1)(i-2k-1)}$. Here we used Hensel's
lemma. One checks that $nk+(n-1)(i-2k-1) \geq \ulcorner (n/2)(i-1)
\urcorner$, and consequently we obtain that the number of elements
of $\mathcal{O}_{(b_1,\ldots,b_n)}$ is a multiple of $q^{\ulcorner
(n/2)(i-1) \urcorner}$.

If $(b_1',\ldots,b_n') \in \mathcal{O}_{(b_1,\ldots,b_n)}$, then
$\mathcal{O}_{(b_1',\ldots,b_n')} =
\mathcal{O}_{(b_1,\ldots,b_n)}$. Consequently, the set $\{
\mathcal{O}_{(b_1,\ldots,b_n)} \mid (b_1,\ldots,b_n) \in
\mathcal{O} \}$ is a partition of $\mathcal{O}$. $\qquad \Box$

\section{The smallest poles of Igusa's zeta function and congruences} \noindent
\textbf{(3.1)} Let $f$ be a rigid $K$-analytic function on $R^n$
which is defined over $R$. Let $S:= \{z/q^i \mid z \in
\mathbb{Z},i \in \mathbb{Z}_{\geq 0} \}$. Using an embedded
resolution of $f$, Igusa proved (see \cite{Igusabook} or
\cite[Section 2.2]{Segers}) that Igusa's $p$-adic zeta function of
$f$ is a rational function in $t$ of the form
\[
Z_f(t) = \frac{A(t)}{\prod_{j \in J} (1-q^{-\nu_j}t^{N_j})},
\]
where $A(t) \in S[t]$, where all $\nu_j$ and $N_j$ are in
$\mathbb{Z}_{>0}$ and where $A(t)$ is not divisible by any of the
$1-q^{-\nu_j}t^{N_j}$.  Remark that the real parts of the poles of
$Z_f(s)$ are the $-\nu_j/N_j$, $j \in J$. Put $l:=\min
\{-\nu_j/N_j \mid j \in J \}$.

It follows from $P(t)=(1-tZ_f(t))/(1-t)$ and $Z_f(t=1)=1$ that the
Poincar\'e series $P(t)=\sum_{i=0}^{\infty} M_i q^{-ni}t^i$ of $f$
can be written as
\[
P(t) = \frac{B(t)}{\prod_{j \in J} (1-q^{-\nu_j}t^{N_j})},
\]
where $B(t) \in S[t]$ is not divisible by any of the
$1-q^{-\nu_j}t^{N_j}$.

In the next paragraphs, we will work in a more general context. By
abuse of notation, we will use the symbols of this particular
situation.

\vspace{0,5cm}

\noindent \textbf{(3.2)} Let $P(t)$ be an \textsl{arbitrary}
rational function in $t$ of the form
\[
P(t) = \frac{B(t)}{\prod_{j \in J} (1-q^{-\nu_j}t^{N_j})},
\]
where $B(t) \in S[t]$, where all $\nu_j$ and $N_j$ are in
$\mathbb{Z}_{>0}$ and where $B(t)$ is not divisible by any of the
$1-q^{-\nu_j}t^{N_j}$. Put $l:=\min \{-\nu_j/N_j \mid j \in J \}$.
Define the numbers $M_i$ by the equality
\[
P(t) = \sum_{i=0}^{\infty} M_i q^{-ni} t^i.
\]
The following proposition is also in \cite{Segers}.

\vspace{0,5cm}

\noindent \textbf{(3.3) Proposition.} \textsl{There exists an
integer $a$ which is independent of $i$ such that $M_i$ is an
integer multiple of $q^{\ulcorner(n+l)i-a\urcorner}$ for all $i
\in \mathbb{Z}_{\geq 0}$.}

\vspace{0,2cm}

\noindent \textsl{Remark.} (i) The statement in the proposition is
obviously equivalent to the following. If $l' \leq l$, then there
exists an integer $a$ which is independent of $i$ such that $M_i$
is an integer multiple of $q^{\ulcorner(n+l')i-a\urcorner}$ for
all $i \in \mathbb{Z}_{\geq 0}$.
\\(ii) Suppose that we are in the situation of (3.1). Then
$n+l>0$, see \cite[Section 3.1.4]{Segersthesis} or (3.4), so that
$(n+l)i-a$ rises linearly as a function of $i$ with a slope
depending on $l$. The statement is trivial if $(n+l)i-a \leq 0$
because the $M_i$ are integers. If $(n+l)i-a>0$, which is the case
for $i$ large enough, it claims that $M_i$ is divisible by
$q^{\ulcorner(n+l)i-a\urcorner}$.

\vspace{0,2cm}

\noindent \textsl{Proof.} We will say that a formal power series
in $t$ has the divisibility property if the coefficient of
$t^i/q^{ni}$ is an integer multiple of
$q^{\ulcorner(n+l)i\urcorner}$ for every $i$.

For $j \in J$, the series
\[
\frac{1}{1-q^{-\nu_j}t^{N_j}}= \sum_{i=0}^{\infty} q^{-i\nu_j}
t^{iN_j} = \sum_{i=0}^{\infty} q^{i(nN_j-\nu_j)}
\frac{t^{iN_j}}{q^{niN_j}}
\]
has the divisibility property because $nN_j-\nu_j$ is an integer
larger than or equal to $N_j(n+l)$. Let $a$ be an integer such
that the polynomial $q^a B(t)$ has the divisibility property.

One can easily check that the product of a finite number of power
series with the divisibility property also has the divisibility
property. This implies that $q^aP(t)$ is a power series with the
divisibility property. Hence $M_i$ is an integer multiple of
$q^{\ulcorner(n+l)i\urcorner -a} = q^{\ulcorner(n+l)i-a\urcorner}$
for all $i$. $\qquad \Box$

\vspace{0,5cm}

\noindent \textbf{(3.4) Proposition.} \textsl{There exist an
integer $a$ which is independent of $i$ and positive integers $R$
and $c$ such that $M_{iR+c}$ is not an integer multiple of
$q^{\ulcorner (n+l)(iR+c)+a \urcorner}$ for $i$ large enough.}

\vspace{0,2cm}

\noindent \textsl{Consequences.} (i) If there exists an integer
$a$ such that $M_i$ is an integer multiple of $q^{\ulcorner
(n+l')i-a \urcorner}$ for all $i \in \mathbb{Z}_{\geq 0}$, then
$l' \leq l$. This is the converse of proposition 3.3. \\ (ii)
Because we saw in the previous section that $M_i$ is an integer
multiple of $q^{\ulcorner(n/2)(i-1)\urcorner}$ if we are in the
situation of (3.1), we obtain already that $Z_f(s)$ has no poles
with real part less than $-n/2$.

\vspace{0,2cm}

\noindent \textsl{Proof.}  Put $J_1=\{j \in J \mid -\nu_j/N_j=l\}$
and $J_2=J \setminus J_1$. Let $N$ be the lowest common multiple
of the $N_j$, $j \in J_1$, and let $\nu$ be the lowest common
multiple of the $\nu_j$, $j \in J_1$. Remark that
$\nu/N=\nu_j/N_j$ for all $j \in J_1$. Let $m$ be the cardinality
of $J_1$. Because $1-q^{-\nu}t^N$ is a multiple of
$1-q^{-\nu_j}t^{N_j}$ for all $j \in J_1$, we can write
\[
P(t) = \frac{D(t)}{(1-q^{-\nu}t^N)^m \prod_{j \in J_2}
(1-q^{-\nu_j}t^{N_j})},
\]
where $D(t) \in S[t]$. Applying decomposition into partial
fractions in $\mathbb{Q}(t)$, we can write
\[
wP(t) = \frac{r_1}{(1-q^{-\nu}t^N)^m} +
\frac{r_2}{(1-q^{-\nu}t^N)^{m-1}} + \cdots +
\frac{r_m}{1-q^{-\nu}t^N} + \frac{E(t)}{\prod_{j \in J_2}
(1-q^{-\nu_j}t^{N_j})},
\]
where $w \in \mathbb{Z}$, where $r_i \in S[t]$ with $\deg(r_i) <
N$ and where $E(t) \in S[t]$.

Let $k \in \mathbb{Z}_{>0}$. Then
\[
\frac{1}{(1-q^{-\nu}t^N)^k} = \sum_{i=0}^{\infty} f_k(i)
q^{-i\nu}t^{iN},
\]
where $f_k : \mathbb{Z}_{\geq 0} \rightarrow \mathbb{Z}_{\geq 0}$
is a polynomial function with rational coefficients of degree
$k-1$. Indeed, because $\sum_{i=0}^{\infty} f_k(i) q^{-i\nu}t^{iN}
= (\sum_{i=0}^{\infty} q^{-i\nu}t^{iN}) (\sum_{i=0}^{\infty}
f_{k-1}(i) q^{-i\nu}t^{iN})$, we obtain that $f_k(n) =
\sum_{i=0}^{n} f_{k-1}(i)$. Consequently, it follows by induction
on $k$ since $\sum_{i=0}^{n} g(i)$ is a polynomial function in $n$
of degree $r$ with rational coefficients if $g$ is such a function
of degree $r-1$.

There exists an integer $d$, an integer $z$ which is not divisible
by $q$ and polynomial functions with integer coefficients $g_b:
\mathbb{Z}_{\geq 0} \rightarrow \mathbb{Z}_{\geq 0}$, $b \in
\{0,1,\ldots,N-1\}$, such that
\begin{eqnarray}
\lefteqn{ \frac{r_1}{(1-q^{-\nu}t^N)^m} +
\frac{r_2}{(1-q^{-\nu}t^N)^{m-1}} + \cdots +
\frac{r_m}{1-q^{-\nu}t^N} } \nonumber \\ & = & \sum_{b=0}^{N-1}
\sum_{i=0}^{\infty} \frac{g_b(i)}{zq^d} q^{-i\nu} t^{iN+b}
\nonumber
\\ & = & \sum_{b=0}^{N-1} \sum_{i=0}^{\infty} \frac{g_b(i)}{z}
q^{(n+l)(iN+b)-d-bl} \frac{t^{iN+b}}{q^{n(iN+b)}}
\end{eqnarray}
and such that $z$ is a divisor of $g_b(i)$ for every $i \in
\mathbb{Z}_{\geq 0}$. Remark that $(n+l)(iN+b)-d-bl$ is an
integer. Because $D(t)$ is not divisible by $(1-q^{-\nu}t^N)^m$,
at least one of the polynomials $g_b$ is different from the zero
polynomial. Fix from now on a $b \in \{0,\ldots,N-1\}$ for which
$g_b$ is different from $0$. Because $g_b$ is a polynomial
function, it has only a finite number of zeros. Let $h$ be a
positive integer which is not a zero of $g_b$. Let $r$ be an
integer such that $g_b(h)$ is not a multiple of $q^r$. Because
$g_b$ is a polynomial with integer coefficients, we obtain for
every $i \in \mathbb{Z}_{\geq 0}$ that $g_b(iq^r+h)$ is not a
multiple of $q^r$. Let $a$ be the greatest integer less than or
equal to $r-d-bl$. Because $(iq^r+h)N+b=iq^rN+hN+b$, we put
$R=q^rN$ and $c=hN+b$. With this notation, we have that the
coefficient of $t^{iR+c}/q^{n(iR+c)}$ in (3) is not an integer
multiple of $q^{(n+l)(iR+c)+r-d-bl}=q^{\ulcorner (n+l)(iR+c)+a
\urcorner}$ for every $i \in \mathbb{Z}_{\geq 0}$.

Now we consider the remaining part
\begin{eqnarray}
\frac{E(t)}{\prod_{j \in J_2} (1-q^{-\nu_j}t^{N_j})}
\end{eqnarray}
of $wP(t)$. We obtain from Proposition 3.3 that there exists an
$l'>l$ and an integer $a'$ such that the coefficient of
$t^i/q^{ni}$ in the power series expansion of (4) is an integer
multiple of $q^{\ulcorner (n+l')i-a' \urcorner}$ for all $i \in
\mathbb{Z}_{\geq 0}$. Consequently, the coefficient of
$t^i/q^{ni}$ is an integer multiple of $q^{\ulcorner (n+l)i+a
\urcorner}$ for $i$ large enough.

Because $wM_{iR+c}$ is the sum of two numbers of which exactly one
is an integer multiple of $q^{\ulcorner (n+l)(iR+c)+a \urcorner}$
for $i$ large enough, we obtain that $wM_{iR+c}$, and thus also
$M_{iR+c}$, is not an integer multiple of $q^{\ulcorner
(n+l)(iR+c)+a \urcorner}$ for $i$ large enough. $\qquad \Box$

\vspace{0,5cm}

\noindent \textbf{(3.5)} Let $\chi$ be a character of $R^{\times}$
with conductor $e$.  Suppose that the image of $\chi$ consists of
the $m$th rooths of unity. Let $\xi= \exp(2 \pi \sqrt{-1}/m)$. The
minimal polynomial of $\xi$ over $\mathbb{Q}$ is the $m$th
cyclotomic polynomial having degree the Euler number $\phi(m)$.
Remark also that $\{ 1,\xi,\ldots,\xi^{\phi(m)-1} \}$ is a basis
of $\mathbb{Q}(\xi)$ as a vector space over $\mathbb{Q}$.

Write
\begin{eqnarray}
Z_{f,\chi}(t) & = & \sum_{i=0}^{\infty} \sum_{u \in
(R/P^e)^{\times}} \chi(u) M_{i+e}(\pi^iu) q^{-n(i+e)}t^i \nonumber
\\ & = & \sum_{0 \leq k < \phi(m)} \left( \sum_{i=0}^{\infty}
\widetilde{M}_{i+e,k} q^{-n(i+e)} t^i \right) \xi^k,
\end{eqnarray}
where $\widetilde{M}_{i+e,k}$ is a linear combination of the
$M_{i+e}(\pi^iu)$, $u \in (R/P^e)^{\times}$, with integer
coefficients because the $m$th cyclotomic polynomial is monic and
has integer coefficients. Since $M_{i+e}(\pi^iu)$ is an integer
multiple of $q^{\ulcorner (n/2)(i+e-1) \urcorner}$, this implies
that $\widetilde{M}_{i+e,k}$ is also an integer multiple of
$q^{\ulcorner (n/2)(i+e-1) \urcorner}$.

On the other hand, using an embedded resolution of $f$, Igusa
proved (see \cite{Igusabook} or \cite[Section 2.2]{Segers}) that
$Z_{f,\chi}(t)$ can be written in the form
\begin{eqnarray}
Z_{f,\chi}(t) & = & \sum_{0 \leq k < \phi(m)}
\frac{A_k(t)}{\prod_{j \in J_k}(1-q^{-\nu_j}t^{N_j})} \xi^k,
\end{eqnarray}
where $A_k(t) \in S[t]$, where every $\nu_j$ and $N_j$ is in
$\mathbb{Z}_{>0}$ and where $A_k(t)$ is not divisible by
$1-q^{-\nu_j}t^{N_j}$ for every $j \in J_k$.

Because $\{ 1,\xi,\ldots,\xi^{\phi(m)-1} \}$ is a basis of
$\mathbb{Q}(\xi)$ as a vector space over $\mathbb{Q}$, we obtain
from (5) and (6) that
\[
\frac{A_k(t)}{\prod_{j \in J_k}(1-q^{-\nu_j}t^{N_j})} =
\sum_{i=0}^{\infty} \widetilde{M}_{i+e,k} q^{-ne} q^{-ni} t^i
\]
for every $k \in \{0,1,\ldots,\phi(m)-1) \}$. The consequence (i)
of (3.4) can be applied to this equality. Because
$\widetilde{M}_{i+e,k}q^{-ne}$ is an integer multiple of
$q^{\ulcorner (n/2)(i+e-1)-ne \urcorner}$, which is equal to
$q^{\ulcorner (n-n/2)i-(n/2)(e+1) \urcorner}$, we obtain that
$Z_{f,\chi}(s)$ has no pole with real part less than $-n/2$.

We obtain the following theorem by using (1.3).

\vspace{0,2cm}

\noindent \textbf{Theorem.} \textsl{Let $n \in \mathbb{Z}_{>1}$.
Let $f$ be a $K$-analytic function on an open and compact subset
of $K^n$. Let $\chi$ be a character of $R^{\times}$. Then we have
that Igusa's $p$-adic zeta function $Z_{f,\chi}(s)$ of $f$ has no
poles with real part less than $-n/2$.}

\vspace{0,2cm}

\noindent \textsl{Remark.} To any $f \in
\mathbb{C}[x_1,\ldots,x_n]$ and $d \in \mathbb{Z}_{>0}$ Denef and
Loeser associate the topological zeta function
$Z_{\mathrm{top},f}^{(d)}(s)$, see \cite{DenefLoeser1} or
\cite{Denefreport}. Because $Z_{\mathrm{top},f}^{(d)}(s)$ is a
limit of Igusa's $p$-adic zeta functions, we obtain that
$Z_{\mathrm{top},f}^{(d)}(s)$ has no poles less than $-n/2$. This
is well known for $n=2$ and we proved this already for $n=3$ in
\cite{SegersVeys}.

\footnotesize{

\noindent \textsc{K.U.Leuven, Departement Wiskunde,
Celestijnenlaan 200B, B-3001 Leuven, Belgium,} \\ \textsl{E-mail:}
dirk.segers@wis.kuleuven.be} \\ \textsl{URL:}
http://wis.kuleuven.be/algebra/segers/segers.htm

\end{document}